\documentclass{aptpub}
\usepackage{epsfig}
\usepackage{graphicx}

\authornames{Bordenave, Foss, Shneer} 
\shorttitle{A Random Multiple Access Protocol} 

\def\ind{{\rm 1\hspace{-0.90ex}1}}

\begin{document}

\title{A Random Multiple Access Protocol with Spatial Interactions} 

\authorone[CNRS and University of Toulouse]{Charles Bordenave} 

\addressone{CNRS UMR 5219,
Institut de Math\'ematiques de Toulouse,
118 route de Narbonne, F-31062 Toulouse, France.
} 

\authortwo[Heriot-Watt University, Edinburgh, UK and Institute of Mathematics, Novosibirsk,
Russia]{Serguei Foss}
\addresstwo{Department of AMS, Heriot-Watt University, Edinburgh, UK, EH14
4AS. 
}

\authorthree[Eindhoven University of Technology and EURANDOM]{Vsevolod Shneer}
\addressthree{EURANDOM, P.O.Box 513-5600 MB, Eindhoven, The Netherlands}

\begin{abstract}
We analyse an ALOHA-type random multiple-access protocol where
users have local interactions. We show that the fluid model of the
system workload satisfies a certain differential equation. We
obtain a sufficient condition for the stability of this
differential equation and deduce from that a sufficient condition
for the stability of the protocol. We discuss the necessary
condition. Further, for the underlying Markov chain, we estimate
the rate of convergence to the stationary distribution. Then we
establish an interesting and unexpected result showing that the
main diagonal is locally unstable if the input rate is
sufficiently small. Finally, we consider two generalisations of
the model.
\end{abstract}



\section{Introduction and stability result}

\subsection{A spatial ALOHA}

We consider a random spatial service system governed by an
ALOHA-type algorithm. More precisely, time is slotted, during a
time slot $n=1,2,..$ a random number $\xi_n$ of users arrive in
the system, and at each slot every user in the system requires
service (transmission) with a certain (transmission) probability
independently of all other users. The sequence $\{ \xi_n\}$ is assumed
to be i.i.d.

The ALOHA multi-access algorithm was first proposed by Abramson
\cite{Abramson}. The slotted scheme was introduced by Roberts
\cite{Roberts}. We consider the latter setting. In the
conventional centralised slotted ALOHA model, there is a single
server. At the beginning of a time slot $n$, the total number
$W_n$ of users in the system is known, and each of them asks for
service (transmission) with probability $\dfrac{1}{W_n}$
independently of all others. If two or more users attempt
transmissions simultaneously, then transmissions collide, the
users remain in the system and try to transmit later. All service
times are equal to $1$, and the server is free at the beginning of
any time slot. It is easy to show that for this system the maximum
throughput is equal to $e^{-1}$. Further, the Markov chain
$\{W_n\}$ is positive recurrent if $\mathbf E \xi_1 < e^{-1}$ and
transient if $\mathbf E \xi_1 > e^{-1}$.

When the information on the numbers $W_n$ of users is unavailable,
various decentralised adaptive algorithms have been introduced and
studied. Algorithms of this type use the information on what has occurred
in the previous time slot: either a conflict or a successful
transmission or no attempted transmissions.
%
More precisely, let $B_n$ be the number of users trying to
transmit during a time slot $n$. In a decentralised algorithm,
only the values of $\min\{B_n,2\}$ are available at time $n+1$.
Assuming that at any time slot $n$ the transmission probability is
the same for all users, Hajek \cite{Hajek} proved that if
$\{\xi_n\}$ are i.i.d. with a finite exponential moment, then
${\bf E} \xi_n < e^{-1} \approx 0.37$ is still necessary and
sufficient for the existence of a stable algorithm. Mikhailov
\cite{Mikhailov} generalised this result by weakening the
exponential moment assumption to the requirement for only the
second moment to exist, while Foss \cite{Foss} generalised it
further by dropping this as well as the independence assumption.
We also refer to Ephremides and Hajek \cite{survey} for a survey
which includes, in particular, results in this direction.


The model described above ignores the network's spatial diversity. In this
paper we present a new model which allows
the possibility for interaction between users to occur
 only when the distance between these users
is small. A limited spatial interaction is a common feature of
wireless networks and is thus very important for practical
consideration.

In this paper we consider only spatial centralised schemes (protocols)
where the total number of users in a neighbourhood is known.
The study of decentralised protocols is a subject of our future research.
Our main result establishes  the stability of the model under
consideration. For that, we use the fluid approximation approach.
We show that the state is repelled by the boundary, and in the interior, the
Euclidean norm turns out to be a Lyapunov function for the fluid model.
Other results of the paper include convergence rates to stationarity and
local (in)stability.

The remainder of the paper is organised as follows. The rest of
this introduction is devoted to the description of the model and
the statement of our main stability result. In Section
\ref{sec:fluid} we prove that fluid limits of the workload in the
system satisfy a certain differential equation.
Section \ref{sec:boundary} is devoted to the analysis of the
behaviour of fluid limits on the boundary of the positive orthant.
In Section \ref{sec_proof} we present a proof of our main
stability result and formulate one of its possible
generalisations. Section \ref{sec:rates} deals with rates of
convergence towards the steady state. Sections \ref{sec:local} and
\ref{behaviour} contain some interesting results on the behaviour
of solutions to the differential equation. Finally, in Section
\ref{sec:extensions} we conclude with some extensions of our
model, which are in a certain sense more applicable to real
systems. In  particular, these extensions include systems where various changes in environment conditions
may result in changes in the radius and/or direction of
interference between the message transmissions).

\subsection{Model description} \label{sec:model}

Let $\mathcal G = (\mathcal V,\mathcal E)$ be an undirected graph
with a finite set of vertices, say $\mathcal V = \{1,...,K\}$. We
suppose that $\mathcal{G}$ is connected. For the graph $\mathcal{G}$ we use
the standard notion of the graph distance. For $i \in \mathcal V$, let $V_i = \{i\} \cup \{ j
\in \mathcal V : \hbox{ such that } (i,j) \in \mathcal E\}$ be the
neighbourhood of a vertex $i$, i.e. the set of vertices at a
maximum distance of $1$ from $i$ in the graph. In the particular case
where all vertices have the same number of neighbours, we denote by $V$
the cardinality of $V_i$, $card(V_i)=V$.

We introduce now a service system with spatial (neighbourhood)
interactions associated with the graph $\mathcal{G}$. We assume
that time is slotted, i.e., arrivals and services may occur only
at times $n=1,2, \ldots$. Suppose that there is a service station
at each vertex of $\mathcal{G}$. The arrival process is denoted by
$A = (A (n) )_{n \in \mathbb N}$, where $A(n) = (A_1(n), ...,
A_K(n)) \in \mathbb N ^ K $, and $A_i(n)$ is the number of users
arriving at time $n$ at a vertex $i$. For $s<t$, denote by $A(s,t)
= \sum_{n= \lceil s \rceil} ^ {\lceil t \rceil -1} A(n)$ the
vector representing the number of users arriving between time instances
$s$ and $t-$. We
suppose that $(A(n))$ is an i.i.d. sequence (while the coordinates
of the vectors may be dependent). We also suppose that $\mathbf E
A_i(n) = \lambda_i
> 0$ for $i = 1, \ldots, K$. If values of $\lambda_i$ do not depend on $i$, then we write
$\lambda_i=\lambda$, $i=1,\ldots, K$.

Let $W (n) = (W_1(n),..., W_K(n)) \in \mathbb{R}_+^K$, where $W_i
(n)$ is the number of users at vertex $i$ at time $n$. At time
$n$, a user at vertex $i$ attempts transmission independently of
the others with probability $1/\sum_{j \in V_i} W_j (n)$. This
user receives service if he is the only user attempting
transmission in $V_i$ at time $n$. We suppose that all service
times are equal to $1$ and that any user leaves the system
immediately upon service completion. Let $N_i (n)$ be the number
of users attempting transmission at time $n$ at vertex $i$. Then
$N_i(n)$ is a binomial random variable with parameters
$\left(W_i(n),\dfrac{1}{\sum_{j \in V_i} W_j (n)}\right)$ and
$\{N_i(n), 1 \leq i \leq K\}$ are independent variables conditioned
on $W(n)$. The sequence $\{W(n)\}$ forms a time-homogeneous
irreducible Markov Chain, which satisfies the recursion

\begin{equation}
\label{eq:evolution} W_i (n) = W_i (n-1) + A_i(n) - \ind (N_i
(n-1) = 1) \prod_{j \in V_i \backslash \{i\}} \ind (N_j (n-1) =
0).
\end{equation}

To explicitly show the dependence of $W(n)$ on the initial
condition $W(0) = x$, we may sometimes write $W^x(n)$.

If $x_i >0$, then the $i$-th component of the drift vector is given by
the expression

\begin{equation}
\mathbf E\bigm[W_i(1) - W_i(0)|W(0) = x\bigm] = \lambda_i - \frac
{x_i}{\sum_{k \in V_i} x_k }\left(1-\frac {1}{\sum_{k \in V_i}
x_k}\right) ^ {x_i -1} \prod_{j \in V_i \backslash \{i\}}
\left(1-\frac {1}{\sum_{k \in V_j} x_k}\right) ^ {x_j},
\end{equation}
and if $x_i = 0$, then $\mathbf E\bigm[W_i(1) - W_i(0)|W(0) =
x\bigm] = \lambda_i$.

We rewrite the expression for the drift vector in the following
way:
$$ \mathbf E\bigm[W(1) - W(0)|W(0) = x \bigm] = \underline \lambda - G(x).$$
Here $\underline{\lambda}$ is
a $K$-dimensional vector with its $i$-th component being equal
to $\lambda_i$ and $G(x)=(G_1(x),\ldots,G_K(x))$
is a function from $\mathbb{R}^K$ to
$\mathbb{R}^K$ defined by
\begin{eqnarray*}
G_i(x) = \begin{cases} \frac {x_i}{\sum_{k \in V_i} x_k
}\left(1-\frac {1}{\sum_{k \in V_i} x_k}\right) ^ {x_i -1}
\prod_{j \in V_i \backslash \{i\}} \left(1-\frac {1}{\sum_{k \in
V_j} x_k}\right) ^ {x_j}, \quad \text{if} \quad x_i > 0, \cr 0,
\quad \text{if} \quad x_i=0.
\end{cases}
\end{eqnarray*}
For $x \in \mathbb R^ K $, we denote $\phi_i (x)=
\frac{x_i}{\sum_{j \in V_i} x_j}$. Let $\phi (x) =
(\phi_1(x),...,\phi_K (x))$.
Note that $G_i$ is bounded from above by $1$ and if $\sum_{k \in V_i} x_k
>0$ then $$ \lim_{t\to +\infty} G_i(tx) = \widetilde G_i(x) =
\phi_i(x) e^{-\sum_{j\in V_i} \phi_j(x)}. $$ In particular,
$\widetilde G_i$ is a homogeneous function of order $0$, i.e.
$\widetilde{G}_i (cx) = \widetilde{G}_i (x)$ for any $c > 0$.

We now comment on the model. In this paper, we mostly consider the
{\it symmetric} case, where $\lambda_i = \lambda$ for all $i = 1,
\ldots, K$ and the graph $\mathcal G$ is $(V-1)$-regular: the
cardinality of $V_i$ is equal to $V$ for all $i$. Notice that, even
in this case, the graph $\mathcal G$ is not necessarily completely
symmetric.



Note also that the system is not monotone. Indeed, $x  \leq y$
(component-wise) does not imply that $W^y(1)$ stochastically
dominates $W^x(1)$. Neither is the system  monotone with
respect to the graph structure: if $G_1$ is embedded into $G_2$,
this does not imply that the workload process built on graph $G_1$
is stochastically dominated by the workload built on graph $G_2$.

We also present a number of generalisations to a non-symmetric
case. In particular, using methods suggested recently in
\cite{Stolyar2}, we formulate in Remark \ref{rem:Stolyar}
sufficient conditions for the stability of the system with
space-inhomogeneous input. Some other generalisations of the model
are proposed in Section \ref{sec:extensions}.

\subsection{Stability result}

We first explain the intuition hidden behind the result.
Consider the symmetric case, $\lambda_i=\lambda$, $card(V_i)=V$,
$i=1,\ldots,K$.

The access protocol favours an equilibrium of the workload in the
network: assume that the workload at node $i$ is much larger than
the workload in its neighbouring nodes, $V_i\setminus \{i\}$.
Then $\phi_i(x)$
is close to $1$, whereas $\phi_j(x)$ is close to $0$ for all the
other
nodes $j$ in $V_i$. Thus the workload at node $j$ in $V_i$ will
tend to get closer to the workload at node $i$. From this balance
mechanism one can guess that the diagonal $\Delta = \{x \in
\mathbb R^K : x_1 = x_2 ... = x_K\}$ is an attractive set.

If the workload vector belongs to the diagonal: $W(0) = c \ind$
where $c \in \mathbb N$, then we obtain:
$$ \mathbf E(W(1) - W(0)|W(0) = c \ind ) = \left(\lambda - \frac 1 V \left(1 - \frac 1 {Vc} \right)^{Vc-1} \right) \ind.$$
Hence, as $c$ tends to infinity, the drift vector converges to
$(\lambda - e^{-1}/V ) \ind$.

So finally, we end up with the conjecture that if $\lambda <
e^{-1}/V$, then the Markov chain $W$ is ergodic. This conjecture
is clearly true for the fully isolated graph and for the complete
graph.

The reasons that led to this conjecture appear to be {\it wrong}
(this will follow from the results of Sections \ref{sec:local} and
\ref{behaviour}, which show that in general
the main diagonal is not attractive).
However, the conjecture itself is true and we can formulate our
main stability result that will be proved in Section
\ref{sec_proof}.

\begin{theorem}
\label{th:main}
Assume the model to be symmetric.

(i) If $\lambda < e^{-1}/V$, then the Markov chain $W$ is positive
recurrent.

(ii) Assume further that $\mathsf P(A_1(1) = A_2(1) =...=A_K(1) =
0) > 0$. Then $W$ is ergodic, i.e. there exists a unique stationary
distribution and the distribution of $W_n$ converges to the
stationary one in the total variation norm.
\end{theorem}

Our proof of (i) is based on the fluid approximation approach. We
will show that all fluid limits satisfy a certain differential
equation and then proceed with the analysis of that equation. The
proof of (ii) is standard: the state ${\bf 0} = (0, 0,...,0)$ is
achievable from any other state and from itself. So, the state
${\bf 0}$ is positive recurrent (due to (i)) and the Markov chain
$W$ is aperiodic. The ergodicity follows.

Our heuristics suggest also that if  $\lambda > e^{-1}/V$ then $W$
is transient. Corollary \ref{cor:eqdiff} in  Section
\ref{sec:local} is a partial result which corroborates this
intuition.

\section{Fluid approximation method}
\label{sec:fluid}

This section deals with a general (not necessarily symmetric)
graph.

\subsection{General properties}

In what follows, we endow $\mathbb R ^K$ with the $L^1$-norm: $|x|
= \sum_{k=1} ^ K |x_k|$. Let $(x^n), n \in \mathbb N,$ be a
sequence in $\mathbb N^K$ such that $\lim_n |x^n| = \infty$.
For $t \in [0,T]$, we define: $$X^n(t) = \frac{W^{x^n}(\lceil
|x^n|t \rceil)}{|x^n|}.$$

To simplify the notation, for $t \in \mathbb R_+$, we set $W(t) =
W(\lceil t \rceil)$.

Let $\mathbb D([0,T],\mathbb R ^K)$ denote the space of c\`adl\`ag
functions from $[0,T]$ to $\mathbb R^K$ endowed with the usual
Skorokhod topology. This means that the distance between the
functions $f_1$ and $f_2$ is defined by
\[
d_T(f_1, f_2) = \inf \sup_{t \in [0,T]} \{|g(t)-t| +
\rho(f_1(g(t)), f_2(t))\},
\]
where $\rho$ is the $L^1$-metric in $R^K$ and the outer infimum is
taken over all monotone continuous functions $g: [0,T] \to [0,T]$
such that $g(0) = 0$ and $g(T) = T$. Denote by $\mathbb D([0,
\infty))$ the space of $R^K$--valued c\`adl\`ag functions on $[0,
\infty)$ with the metric
\[
d(f_1,f_2) = \sum_{T=1}^\infty 2^{-T} \frac{d_T(f_1,
f_2)}{1+d_T(f_1,f_2)}.
\]
Note that $X^n \in \mathbb D([0,T],\mathbb R ^K)$, for all $n$.

\begin{lemma} \label{le:tight}

(i) For any sequence $x^n$ such that $|x^n| \to \infty$, the
family $ \{ (X^n), n \in \mathbb N\}$ almost surely
has a compact closure in
the Skorokhod topology, and any accumulation point $z$ of $\{X_n\}$
is almost surely continuous.

(ii) Function $z$ is Lipschitz with the constant
$\max\{\sum_{i=1}^K \lambda_i, K\}$.
\end{lemma}

\begin{proof}

(i) One can easily obtain a proof of this statement along the
lines of the proof of \cite{Dai}, Theorem 4.1 or \cite{Stolyar},
Theorem 7.1. Formally, the proofs of the mentioned theorems are
given for multi-class networks. However, as pointed out in
\cite{FK}, the tightness of such families holds under weaker
conditions (see \cite{FK}, Assumption 2.19).

(ii) Since $G_i$ are bounded by $1$ from above,

\begin{eqnarray*}
|X^n (t) - X_i^n (s)| & \leq & \max \left\{ \frac{|A(s|x^n|,
t|x^n|)|}{|x^n|}, \frac{K |x^n| (t-s+1/|x^n|)}{|x^n|}\right\}
\\ & \leq &  \max\left\{ \frac{1}{ |x^n|}
\sum_{k =\lceil |x^n| s \rceil }^{\lfloor |x^n| t \rfloor } B_k  ,
K(t-s+1/|x^n|)\right\},
\end{eqnarray*}
where $B_k$ is the total number of arrivals at time $k$. Sequence
$\{B_k\}_{k \in \mathbb N}$ consists of i.i.d. random variables
with $\mathbf E B_k = K \overline \lambda$. By the law of large
numbers, the result now follows if we let $n \to \infty$.

\end{proof}

\begin{definition}
Any accumulation (in the Skorokhod topology) point $z$ of the
sequence $X^n$ is a {\it fluid limit}. The collection of all fluid
limits is called the {\it fluid model}.
\end{definition}

It follows from the definition of $X^n$ and $z$ that $|z(0)|=1$
and that $z_i(t) \ge 0$ for all $i=1, \ldots, K$ and for all $t$.

\begin{corollary}
The sample-path trajectories of fluid limits are self-similar.
More precisely, for any fluid limit $z$ and for any $u > 0$ such
that $\mathsf P(|z(u)| > 0) > 0$, the random process $\{
\widetilde{z}(t), t \ge 0\}$ with conditional distribution
\[
\mathsf P(\widetilde{z} (t) \in \cdotp) = \mathsf
P\left(\frac{z(u+t)}{|z(u)|} \in \cdotp \big| z(u)\right)
\]
is also a fluid limit on the event $|z(u)| > 0$.
\end{corollary}

This result may be obtained along the lines of the proof of
Stolyar \cite{Stolyar}, Lemma 6.1. However, the remark given in
the proof of Lemma \ref{le:tight} (i) also applies here.

\begin{definition}
We say that the fluid model is {\em stable} if there exists a
deterministic time $t_0$ and $\varepsilon \geq 0$, such that $|z
(t)| \leq \varepsilon$ a.s. for $t \geq t_0$, for all fluid limits
$z$. Due to the self-similarity of fluid limits, this is
equivalent to saying that there exists a deterministic
time $t$ such that $|z(t)|=0$ a.s., for all fluid limits $z$.
\end{definition}


Definition 2 of fluid stability has become standard and appears in
most papers dealing with the fluid approximation method.

\subsection{Fluid model criterion for stability}

In this subsection we formulate a stability theorem for fluid
limits which will imply Theorem \ref{th:main} (i).

\begin{lemma}
\label{th:fluid2chain}
\begin{itemize}
If the fluid model is stable then $W$ is positive recurrent.
\end{itemize}
\end{lemma}

\begin{proof} One can again obtain a proof of this assertion by
following the lines of the proofs of Dai \cite{Dai} or Stolyar
\cite{Stolyar} which are given for multi-class networks.


\end{proof}

By Lemma \ref{th:fluid2chain}, Theorem \ref{th:main} (i) will
follow from the next statement.

\begin{theorem}
\label{th:fluidstable} If $\lambda < e^{-1}/V$, then the fluid
model is stable.
\end{theorem}

Our proof of Theorem \ref{th:fluidstable} is given in Section
\ref{sec_proof}.

\subsection{Fluid limit evolution equation}

In what follows we write
 \begin{equation}\label{whatfollows}
\varphi_i(t) = \phi_i (z(t)) =
\dfrac{z_i(t)}{\sum\limits_{j \in V_i} z_j(t)}.
\end{equation}

\begin{theorem}
\label{proofderivative} Take any fluid limit $z$ and fix $t \ge
0$. Assume that $\sum_{j \in V_i} z_j(t) > 0$ for all $i$. Then
$z_i(t)$ is differentiable at point $t$ if $t>0$, and has a right
derivative at point $t$ if $t=0$. Moreover,
\begin{equation}
\label{eq:diffeq} z_i^{'}(t) = \lambda_i - \varphi_i(t) e^{-
\sum\limits_{j \in V_i} \varphi_j(t)} = \lambda_i - \widetilde G_i
(z_i(t)),
\end{equation}
where $z_i(t)$ is the (right) derivative.
\end{theorem}

Under the assumptions of the Theorem, this differential equation
admits a unique solution. 


Fluid limits with an initial condition on the boundary ($\sum_{j
\in V_i} z_j(0) = 0$ for some $i$) are analysed in Section
\ref{sec:boundary}.

\begin{proof}
[Proof of Theorem \ref{proofderivative}]

Two cases are possible: either $z_i(t) > 0$, or $z_i(t)=0$ and
$\sum_{j \in V_i} z_j(t) > 0$. We treat these two cases
separately.

(i) Suppose first that $z_i (t) >0$. To treat this case, we need
the following technical result.

\begin{lemma} \label{le:FF} There exists $C>0$ such that
$|G_i(x) - \widetilde G_i(x)| \leq \min(1,C /x_i)$ if $x_i \ge 2$.
\end{lemma}

\begin{proof}
[Proof of Lemma \ref{le:FF}]

Using that $|e^{-y_1} - e^{y_2}| \le |y_1-y_2|$ for all $y_1, y_2
\ge 0$, we obtain
\begin{eqnarray}
|G_i(x) - \widetilde G_i(x)| & \leq & \left|\ln\left(1-\frac
{1}{\sum_{k \in V_i} x_k}\right)\right| \nonumber \\ & + &
\left|\sum_{j \in V_i} \left(x_j \ln\left(1-\frac {1}{\sum_{k \in
V_j} x_k}\right)+\frac {x_j}{\sum_{k \in V_j} x_k}\right)\right|.
\end{eqnarray}
For every $j$, denote $y_j = \dfrac{1}{\sum_{k \in V_j} x_k}$.
From the inequality $|\ln(1-y) + y| \leq \dfrac{y^2 }{2(1-y)^2}$
for $y \in (0,1)$  we obtain that
\begin{eqnarray*}
|G_i(x) - \widetilde G_i(x)| & \le & y_i + \frac{y_i^2}{2
(1-y_i)^2} + \sum_{j \in V_i}  \frac{x_j y_j^2}{2 (1-y_j)^2}.
\end{eqnarray*}
The statement of the lemma now follows from the inequalities
\[
y_j \le 1/x_i, \quad x_j y_j \le 1 \quad \text{and} \quad y_j \le
1/2
\]
for all $j \in V_i$.

\end{proof}

Assume now that $t=0$ (the result for an arbitrary $t$ follows
from the self-similarity of fluid limits). Let $x = z(0)$. Suppose
that $s < x_i$. Let $k \leq |x^n|s$, then $W_i^{x^n} (k)\geq x^n_i
- k \geq |x^n|(x^n_i/|x^n| -s)$. Hence, $W_i^{x^n} (k) \ge 2$ for
$k \leq |x^n|s$, for large enough $n$.

We need to show that $\lim_{s \to 0} \dfrac{z_i(s) - z_i(0)}{s} =
\lambda_i - \widetilde G_i(z(0))$. Consider the decomposition
\begin{eqnarray}
X_i^n (s) -X_i^n (0) & = & \frac{1}{|x^n|} \sum_{k=0} ^ {\lfloor
|x^n|s \rfloor-1} \left(W_i^{x^n} (k+1) - W_i^{x^n}(k)\right)
\nonumber
\\ & = & \frac{1}{|x^n|}  \sum_{k=0} ^ {\lfloor |x^n|s \rfloor-1}
\mathbf E\biggl[W_i^{x^n} (k+1) - W_i^{x^n}(k)| W^{x^n}(k)\biggr]
\nonumber
\\ && \qquad + \frac{1}{|x^n|}  \sum_{k=0} ^ {\lfloor |x^n|s
\rfloor-1} \left(W_i^{x^n} (k+1) - \mathbf E[W_i^{x^n} (k+1)|
W_i^{x^n}(k)]\right) \nonumber
\\ & = & \frac{1}{|x^n|}  \sum_{k=0} ^ {\lfloor |x^n|s \rfloor-1}
\left(\lambda_i - G_i(W^{x^n}(k))\right) + \frac{1}{|x^n|}
\sum_{k=1} ^ {\lfloor |x^n|s \rfloor} D^n _k, \label{eq:disdiff}
\end{eqnarray}
where
\[
D^n_k = W_i^{x^n}(k) - \mathbf E
\left(W_i^{x^n}(k)|W^{x^n}(k-1)\right) = A_i(k) - \lambda_i +
q_i(k) - \mathbf E \left(q_i(k) |W^{x^n}(k-1)\right)
\]
with $q_i(k) = I (N_i (k-1) = 1) \prod_{j \in V_i \backslash
\{i\}} I (N_j (k-1) = 0)$. We have $\dfrac{1}{|x^n|} \sum_{k=1}^
{\lfloor |x^n|s \rfloor} (A_i(k) - \lambda_i) \to 0$ a.s. as $n
\to \infty$. So we can apply Theorem VII.3 of Feller \cite{feller}
(with $b_k = 1/k$) to deduce that
\begin{equation}
\label{eq:Dnk} \frac{1}{|x^n|} \sum_{k=1}^ {\lfloor |x^n|s
\rfloor} \left(q_i(k) - \mathbf E\left(q_i(k)
|W^{x^n}(k-1)\right)\right) \to 0 \quad \text{a.s.}
\end{equation}
as $n \to \infty$.

It remains to find the limit of the first term in the RHS of
equation (\ref{eq:disdiff}). Decompose this term as follows:

$$\frac{1}{|x^n|}  \sum_{k=0} ^ {\lfloor |x^n|s \rfloor-1}
\left(\lambda_i
 - G_i(W^{x^n}(k))\right) = \frac{1}{|x^n|}  \sum_{k=0} ^ {\lfloor
 |x^n|s
\rfloor-1} \left(\lambda_i - \widetilde G_i\left(X^n
\left(\frac{k}{|x^n|}\right)\right)\right) + \varepsilon(s,n),$$
where by Lemma \ref{le:FF}
\[
|\epsilon(s,n)| \le C \frac{1}{|x^n|} \sum_{k=0}^{\lfloor |x^n|s
\rfloor-1} \frac{1}{W_i^{x^n}(k)} \le C \frac{1}{|x^n|}
\sum_{k=0}^{\lfloor |x^n|s \rfloor-1} \frac{1}{x_i^n-k} \to 0
\]
as $n \to \infty$ uniformly in $s \leq x_i$. Further, from the
uniform convergence of $X^n$ to $z$ and the continuity of
$\widetilde G$ we deduce that
\[
\frac{z_i(s)-z_i(0)}{s} = \lambda_i - \lim_{n \to \infty}
\frac{\sum\limits_{k=0}^{\lfloor |x^n|s \rfloor-1}
\widetilde{G}_i\left(z\left(\frac{k}{|x^n|}\right)\right)}{|x^n|s}.
\]
Since $\frac{1}{|x^n|} \sum_{k=1}^{\lfloor |x^n|s \rfloor}
\widetilde G_i (z(\frac{k-1}{|x^n|}))$ is a Riemann sum of a
continuous bounded function, it converges to $\int_0^s \widetilde
G_i (z(u))du$, so
\begin{equation}
\label{eq:Eeqdiff} \lim_{s \to 0} \frac{z_i(s) - z_i(0)}{s} =
\lambda_i - \lim_{s \to 0} \frac{\int_0^s \widetilde
G_i(z(u))du}{s} = \lambda_i - \widetilde G_i(z(0)).
\end{equation}

(ii) Now consider the second case, $z_i (0) = 0$ and $\sum_{j \in
V_i} z_j(0)
> 0$. Notice that $\widetilde G_i(z_i (0)) = 0$. In view of equations
(\ref{eq:disdiff}) and (\ref{eq:Dnk}) it suffices to show that
\begin{equation}
\label{eq:0} \lim_{s \to 0^+} \lim_{n \to \infty} \frac{1}{|x^n|s}
\sum_{k=0} ^ {\lfloor |x^n|s \rfloor-1}
 G_i(W^{x^n}(k))  = 0.
\end{equation}

By the assumption, there exists $j \in V_i$ such that $z_j(0) =
\lim_{n
  \to \infty} x_j ^ n / |x^n| > \alpha > 0$. Let $\varepsilon>0$, then there exists
$n_0$ such that for all $n \geq n_0$,  $x_j ^ n / |x^n| > \alpha$
and $x_i ^ n / |x^n| < \varepsilon$.

Fix $0< s < \alpha$ and $\varepsilon < \alpha$. Then for $n$ large
enough, $W_i^{x^n}(k) \leq
 \varepsilon|x^n|  + A_i (0,k)$, $W_j^{x^n}(k) \geq \alpha |x^n| - k$ and
\begin{eqnarray*}
 G_i(W^{x^n}(k)) & \leq &  \frac{W_i^{x^n}(k)}{W_i^{x^n}(k) +
   W_j^{x^n}(k)} \\
& \leq &  \frac{ \varepsilon|x^n|  + A_i (0,k )}{(\alpha   +
\varepsilon) |x^n| - k}
\end{eqnarray*}

By the strong law of large numbers, $\lim_{t \to +\infty} A_i
(0,t) / t = \lambda_i$ a.s. Let $\widetilde \lambda > \lambda_i$.
We may choose $|x|$ and $n$ large enough and then $k_0$ such that
for $k_0 \leq k \leq s |x^n|$,
$$
G_i(W^{x^n}(k)) \leq \frac{\varepsilon |x^n| + \widetilde \lambda
k }{(\alpha
   + \varepsilon)  |x^n| - k },
$$
and
\begin{eqnarray*}
\frac{1}{|x^n|}  \sum_{k=0} ^ {\lfloor |x^n|s \rfloor-1}
 G_i(W^{x^n}(k)) & \leq &  \frac{k_0}{ |x^n|} + \frac{\varepsilon}{\alpha +
   \varepsilon- s}  + \frac{1}{|x^n|} \sum_{k=0} ^ {\lfloor |x^n|s \rfloor-1} \frac{k\widetilde \lambda}{\alpha|x^n|
  - k} \\
\end{eqnarray*}
Direct computations show that
$$\lim_{n \to \infty} \frac{1}{|x^n|} \sum_{k=0} ^ {\lfloor |x^n|s \rfloor-1} \frac{k\widetilde \lambda}{\alpha|x^n|
  - k} = -\widetilde \lambda (s+\alpha \ln (1- \frac{s}{\alpha})).
$$
Then
$$
\limsup_{n} \frac{1}{|x^n|}  \sum_{k=0} ^ {\lfloor |x^n|s
\rfloor-1}
 G_i(W^{x^n}(k)) \leq  \frac{\varepsilon}{\alpha +
   \varepsilon- s} - \widetilde \lambda (s + \alpha \ln (1-
   \frac{s}{\alpha}))\quad \text{a.s.}
$$
Since the last inequality holds for all $\varepsilon> 0$ and
$\widetilde \lambda
> \lambda_i $, we have
$$
\limsup_{n} \frac{1}{|x^n|}  \sum_{k=0} ^ {\lfloor |x^n|s
\rfloor-1}
 G_i(W^{x^n}(k)) \leq  -\lambda_i(s+ \alpha \ln (1- \frac{s}{\alpha})).
$$
It then follows immediately that
$$
\lim_{s \to 0^+} \limsup_{n} \frac{1}{|x^n|s}  \sum_{k=0} ^
{\lfloor |x^n|s \rfloor-1}
 G_i(W^{x^n}(k)) = 0.$$
The proof of Theorem 3 is now complete.
\end{proof}

We also need a further result that may be deduced from Theorem
\ref{proofderivative}. Let $H = \{x \in \mathbb R^K: x_i > 0
\quad \text{for all}
\quad i=1, \ldots, K\}$
be the interior of the positive orthant. 

\begin{lemma}
\label{le:trajH} Assume that $z(0) \in H$, then either
\begin{enumerate}
\item[(i)] there exists $c$ such that $z(c) = 0$ and $z(t) \in H$
for all $t \in (0,c)$ or \item[(ii)] $z(t)$ remains in $H$ for all
$t>0$.
\end{enumerate}
\end{lemma}

\begin{proof}
[Proof of Lemma \ref{le:trajH}]

Restricted on the open set $H$, the RHS of equation
(\ref{eq:diffeq}) is a continuous function. Therefore, the
 solutions of equation (\ref{eq:diffeq}) are locally
uniquely defined as long as $z(t)$ remains in $H$. Now, suppose,
to the contrary, that $t \mapsto z(t)$ leaves $H$ at time $c$ at a
point $y = \lim_{t \to c-} z(t) \in \partial H \backslash \{0\}$.

Let $a_i = \limsup_{t \to c^-} \phi(z(t))$, $a_i \in [0,1]$. Since
$y \neq 0$, there exist $i_1$ and $i_2$ such that $y_{i_1} = 0$
and $y_{i_2} >0$. The connectivity of $G$ implies that there
exists $k$ such that $y_k = 0$ and $\sum_{j \in V_k} y_k >0$ (for
that, consider any path from $i_1$ to $i_2$). Hence, $a_k = 0$ and
$\lim_{t \to c^-} F_k (\phi(z(t))) = \lambda > 0$, this implies
that $t \mapsto x_k(t)$ increases in a left neighbourhood of $c$,
which contradicts $y_k = \lim_{t \to c-} z(t) = 0$.
\end{proof}

Lemma \ref{le:trajH} implies that for an initial condition in $H$
the fluid limit $z(t)$ remains in $H$ or finally reaches $0$ at
time $c$. 
It also implies that if $z(0) = \lim_n x^n / |x^n| \in H$ then the
fluid limit is deterministic.

\section{Properties of the fluid limit on the boundary}
\label{sec:boundary}

In this section we work with the general case (we no longer make the symmetry assumptions).

\begin{conjecture}
We conjecture that all coordinates of
any fluid limit $z$ have right derivatives at
point $0$ (even if there
exists $i$ such that $x_j = z_j(0) = 0$ for all $j \in V_i$). We
also conjecture that the right derivative $z^{'}(0)$ of a fluid
limit $\{z(t), t\ge 0\}$ does not depend on the sequence $x^n$ and
only depends on $x = \lim_n x^n/|x^n|$. If this is true, then all fluid
limits are deterministic functions.
\end{conjecture}

In this Section, we prove a weaker statement which implies that
the boundary of the positive orthant does not play any role in
determining the stability conditions for the fluid model. Denote
$$\tau_h = \inf\{t \ge 0: |z(t)| < h\}.$$ Denote also $\lambda_* =
\min\{\lambda_1, \ldots, \lambda_K\}> 0$.
 Since $|z(0)| = 1$, $\max_i z_i(0) \ge 1/K$.
The inequality $z_i^{'}(t) \ge \lambda_* - 1$ for all $i$ and $t$
also implies that
\begin{equation} \label{T_2}
\tau_{1-\varepsilon} \ge \frac{\varepsilon}{K (1-\lambda_*)}.
\end{equation}

\begin{theorem}
\label{achieve} There exist positive constants $b$ and
$\varepsilon_0$ such that, for any $\varepsilon \in
(0,\varepsilon_0)$, $\min_i z_i(t) \ge b \varepsilon $ for any $t
\in [c \varepsilon, \tau_{1-\varepsilon})$ where $c =
1/K(1-\lambda_*)$.
\end{theorem}

From the self-similar property of fluid limits, the following corollary is immediate.
\begin{corollary}
\label{cor:liminf} For any $ h >0$, $z_i(t) > 0$ for all $0 < t < \tau_{h}$
and all $i$.
\end{corollary}

Lemma \ref{le:trajH} and Corollary \ref{cor:liminf} imply the
following.
\begin{corollary}
\label{le:trajH1} Assume that $|z(0)| > 0$, then either
\begin{enumerate}
\item[(i)] there exists $c$ such that $z(c) = 0$ and $z(t)$ remains
in $H$ for all $t \in (0,c)$ or
\item[(ii)] $z(t)$ remains in $H$ for
all $t>0$.
\end{enumerate}
\end{corollary}

We do not present proofs of these statements here as they are
rather obvious.


The rest of this Section is devoted to the proof of Theorem
\ref{achieve}. We begin with two technical lemmas.

\begin{lemma} \label{derivative}
There exist positive constants $K_1 > 1$ and $K_2$ such that, for
any fluid limit $z$, if $z_i(t) > K_1 z_j(t)$ for two neighbouring
nodes $i$ and $j$, then $z_j^{'}(t) > K_2$.
\end{lemma}

\begin{proof}
[Proof of Lemma \ref{derivative}]

The existence of $z_j^{'}(t)$ is guaranteed by Theorem
\ref{proofderivative}. Since $\sum_{k
  \in V_j}  z_k (t) > z_i(t) > 0$, we have
\[
z_j^{'}(t) > \lambda_* - \frac{z_j(t)}{\sum_{k \in V_j} z_k (t)}
\ge \lambda_* - \frac{z_j(t)}{z_i(t) + z_j(t)} > \lambda_* -
\frac{1}{1+K_1}
\]
and we may take $K_1 = 2/\lambda_* - 1$ and $K_2 = \lambda_*/2$.

\end{proof}

\begin{lemma} \label{exit}
There exist constants $C_1 \ge C_2 > 0$ such that for any $h > 0$
one can choose $h_1 > 0$ such that if $|z(0)| \ge h_1$ and $\min_i
z_i(0) \ge C_1 h$ then $\min_i z_i(t) \ge C_2 h_1$ for all $t \le
\tau_{h}$.
\end{lemma}

\begin{proof}
[Proof of Lemma \ref{exit}]
Denote by $D$ the
maximum graph distance in $\mathcal{G}$ (the diameter of
$\mathcal{G}$).  Put $C_1 =
\dfrac{1}{K K_1^{D+1}}$ and $C_2 = \dfrac{C_1}{K_1^{D-1}}$. We may
prove Lemma \ref{exit} for $h = 1$. The result for an arbitrary $h$
follows from the self-similarity of fluid limits.

It is sufficient to show that for any $t < \tau_{1}$ if $\min_i
z_i(t) \ge C_1$ then there exists $0 < s < \infty$ such that
\begin{equation} \label{stat1}
\min_i z_i(t+s) \ge C_1
\end{equation}
and
\begin{equation} \label{stat2}
\min_i z_i(u) \ge C_2 \quad \text{for all} \quad t \le u \le t+s.
\end{equation}
Indeed, assume that (\ref{stat1})-(\ref{stat2}) hold and Lemma
\ref{exit} is not valid. Then there exists $t \le \tau_{1}$ such
that $\min_i z_i(t) < C_2$. It then follows from the continuity of
fluid limits that there is a last moment $v < t$ when
$\min_i z_i(v) \ge C_1$. However, it follows from
(\ref{stat1})-(\ref{stat2}) that there exists $s
> 0$ such that $\min_i z_i(v+s) \ge C_1$ and $\min_i z_i(u)
\ge C_2$ for all $v \le u \le v+s$. Clearly, $v+s < t$, which
contradicts our assumption that $v$ is a last moment before $t$
when $\min_i z_i(v) \ge C_1$.

Now let $t$ be any time  such that $t < \tau_{1}$ and $\min_i
z_i(t) \ge C_1$. Note that $\max_i z_i(t) \ge 1/N = C_1 K_1^{D+1}$
since $t < \tau_{1}$. To simplify the notation, assume that
$z_1(t) = \max_i z_i(t)$. Let $T$ be such that $z_1(t+u) \ge C_1
K_1^D$ for all $0 \le u \le T$. Again note that $z_i^{'}(u) \ge
\lambda_* - 1$ for all $i$ and $u$. This implies that
\begin{equation} \label{T}
T \ge \frac{C_1(K_1^{D+1} - K_1^D)}{1-\lambda_*} = \frac{C_1 K_1^D
(K_1 - 1)}{1-\lambda_*}.
\end{equation}
Let $d$ be the maximum distance in $\mathcal{G}$ from node $1$.
Clearly, $d \le D$. For $j=1, \ldots, d$, denote by $A_j$ the set
of nodes at distance $j$ from node $1$.

We show that there exists $0 < s < T$ such that (\ref{stat1}) and
(\ref{stat2}) hold. First, we show that $\min z_i(u) \ge C_2$ for
all $t \le u \le t+T$. Note that $z_i(u) \ge C_1$ for all $i \in
A_1$ and $t \le u < t+T$. Indeed, assume that there exist $i \in
A_1$ and $t \le u < t+T$ such that $z_i(u) < C_1$. Then, by
continuity, there is a last moment $t \le u_1 < u$ such that
$z_i(u_1) \ge C_1$. Lemma \ref{derivative} implies that
$z_i^{'}(u_1) \ge K_2 > 0$ and hence, there exists time $u_2 >
u_1$ such that $z_i(u_2) \ge C_1$, but this contradicts our
assumption on $u_1$. Using induction and following the same
arguments, we can show that $z_i(u) \ge C_1/K_1^{j-1}$ for all $i
\in A_j$ and $t \le u \le t+T$ for any $j=1, \ldots, d$. Hence,
$\min_i z_i(u) \ge C_1/K_1^{d-1} \ge C_1/K_1^{D-1} = C_2$ for all
$t \le u \le t+T$.

Let us now show that there exists $0 < s < T$ such that
(\ref{stat1}) holds. For every $j=1, \ldots, d$, denote by $t_j$
the time needed to achieve the level $C_1 K_1^{d-j}$ starting from
the level $C_1/ K_1^{j-1}$ and moving with speed $K_2$. Clearly,
$t_j = \dfrac{C_1 (K_1^{d-1}-1)}{K_2 K_1^{j-1}}$. Note that
(\ref{stat1}) and (\ref{stat2}) hold with $s = \sum\limits_{j=1}^d
t_j$ if $T \ge \sum\limits_{j=1}^d t_j$. Indeed, $\min_{j \in A_1}
z_j$ will achieve the level $C_1 K_1^{d-1}$ not later than at time
$t+t_1$ and will not become smaller than this level before time
$t+T$, since all nodes in $A_1$ are neighbours of node $1$ and
$z_1(u) \ge K_1^D$ for all $t \le u \le t+T$. Note also that
$\min_{j \in A_2} z_j$ will become greater than $C_1 K_1^{d-2}$
not later than at time $t+t_1+t_2$ since it cannot become smaller
than $C_1/K_1$ before time $t+t_1$, and after this time it is
either greater than $C_1 K_1^{d-2}$ or grows with a speed of at
least $K_2$ (this follows from Lemma \ref{derivative} and the fact
that any node in $A_2$ has a neighbour in $A_1$). We can continue
these arguments to prove that $\min_{j \in A_d} z_j$ will become
greater than $C_1$ not later than at time $t+\sum\limits_{i=1}^d
t_i$ if $T \ge \sum\limits_{i=1}^d t_i$.

Note that
\begin{eqnarray} \label{sum}
\sum\limits_{i=1}^d t_i & = & \frac{C_1 (K_1^{d-1} -
1)(1+K_1+\ldots+K_1^{d-1})}{K_2 K_1^{d-1}} = \frac{C_1
(K_1^{d-1}-1)(K_1^d-1)}{K_2 K_1^{d-1} (K_1-1)} \cr & \le &
\frac{C_1 (K_1^d-1)}{K_2 (K_1-1)} \le \frac{C_1 (K_1^D-1)}{K_2
(K_1-1)}.
\end{eqnarray}

If we take $K_2 = \lambda_*/2$ and $K_1 = 2/\lambda_* - 1$ then
$(1-\lambda_*)/K_2 = K_1 - 1$. Note also that in this case $K_1
\ge 2$. It now follows from (\ref{T}) and (\ref{sum}) that $T \ge
\sum\limits_{i=1}^d t_i$.

\end{proof}

One can see from the proof of Lemma \ref{exit} that the following
(stronger) result holds.

\begin{lemma} \label{arbith}
For any $h_1 > 0$ there exists $\widehat{h_2} > 0$ such that for
any $h_2 \le \widehat{h_2}$ there exists $0 < h_3 \le h_2$ such
that if $|z(0)| \ge h_1$ and $\min_i z_i(0) \ge h_2$ then $\min_i
z_i(t) \ge h_3$ for all $t \le \tau_{h_1}$.
\end{lemma}

\begin{remark} \label{value}
Lemma \ref{arbith} is valid with $\widehat{h_2} = \dfrac{h_1}{K
K_1^{D+1}}$.
\end{remark}

\begin{proof}
[Proof of Theorem \ref{achieve}]

The proof of Theorem \ref{achieve} is similar to that of Lemma
\ref{exit}.  Recall that $D$ is the maximum graph distance of $\mathcal{G}$. Take $\varepsilon_0$ such that $\dfrac{K_2 (K_1 - 1)
\varepsilon}{(K_1^D - 1)} \le \dfrac{1-\varepsilon}{K K_1^{D+1}}$
for all $\varepsilon \le \varepsilon_0$ and let $a = \dfrac{K_2
(K_1 - 1)}{(K_1^D - 1)}$. Then $a \varepsilon \le
\dfrac{1-\varepsilon}{K K_1^{D+1}}$, and, in view of Lemma
\ref{arbith} and Remark \ref{value}, it is enough to prove that
$\min_i z_i(c \varepsilon) \ge a \varepsilon$.

Note that $\max_i
z_i(0) \ge 1/K$. Assume that $z_1(0) = \max_i z_i(0)$.
Let $T$ be such that $z_1(u) \ge a \varepsilon K_1^D$ for all $0
\le u \le T$. Then $z_i^{'}(t) \ge \lambda_* - 1$
implies that
\begin{equation} \label{T_1}
T \ge \frac{1/K - a \varepsilon K_1^D}{1-\lambda_*} = \frac{1-K a
\varepsilon K_1^D}{K(1-\lambda_*)}.
\end{equation}

Again let $d$ be the maximum distance in $\mathcal{G}$ from node
$1$. 
For $j=1, \ldots, d$, denote by $A_j$ the
set of nodes at distance $j$ from node $1$. For every $j=1,
\ldots, d$, denote by $t_j$ the time needed to achieve the level
$a \varepsilon K_1^{d-j}$ starting from the level $0$ and moving
with speed $K_2$. Clearly, $t_j = \dfrac{a \varepsilon
K_1^{d-j}}{K_2}$.
Denote $T_1 = \sum\limits_{j=1}^D t_j$. Note that
\begin{equation} \label{sum1}
T_1 = \frac{a \varepsilon (K_1^D - 1)}{K_2 (K_1 - 1)} =
\frac{\varepsilon}{K (1- \lambda_*)} = c \varepsilon.
\end{equation}

Following the same arguments as in the proof of Lemma \ref{exit},
we can show that $\min_i z_i(c \varepsilon) = \min_i z_i(T_1) \ge
a \varepsilon$ if $T_1 \le T$.

It remains to prove that $T_1 \le T$. This follows from
(\ref{T_1}), (\ref{sum1}) and the inequality $a \varepsilon \le
\dfrac{1-\varepsilon}{K K_1^{D+1}}$.
\end{proof}

\begin{remark}
Denote by $\nu(s, h, b) = \inf\{t\ge s: |z(t)| < h \quad \text{or}
\quad \min_i z_i(t) < b\}$ the time (after moment $s$) of the
first exit from the set $\{|z| \ge h\} \cap \{\min_i z_i \ge b\}$
after time $s$. Theorem \ref{achieve} implies that there exist $b
> 0$ and $s \ge 0$ such that $\tau_{1-\varepsilon} = \nu(s \varepsilon , 1-\varepsilon, b \varepsilon)$ for
any initial condition $z(0)$.
\end{remark}

\section{Proof of Theorem \ref{th:fluidstable}} \label{sec_proof}

In this Section we first present a proof of our main stability
result and then formulate its generalisation. Recall that Theorem
1 follows from Theorem \ref{th:fluidstable} and Lemma
\ref{th:fluid2chain}.

Recall also that here we deal with the symmetric case. We need to
prove that there exists a deterministic time $t_0$ such that for
all fluid limits, $z (t) =0 $ for $t \geq t_0$ a.s.

\begin{lemma} \label{le:lyapunov}
If $z_i(t) > 0$ for all $i=1, \ldots, K$, then
\[
\left(\sum\limits_i z_i^2(t)\right)' \le \left(\lambda -
\frac{e^{-1}}{V}\right) \sum_i z_i(t)
\]
and hence, if $\lambda < \frac{e^{-1}}{V}$,
\[
\left(\sum\limits_i z_i^2(t)\right)' \le - \varepsilon \sum_i
z_i(t)
\]
for some $\varepsilon > 0$.
\end{lemma}

\begin{proof}
[Proof of Lemma \ref{le:lyapunov}.]

Clearly, it is sufficient to prove the inequality
\begin{equation} \label{inequality}
\frac{\sum\limits_i z_i \varphi_i \exp\left\{-\sum\limits_{j \in
V_i} \varphi_j\right\}}{\sum\limits_k z_k} \ge \frac{e^{-1}}{V}
\end{equation}
where we slightly abuse the notation from
(\ref{whatfollows}) by writing $z_i$ instead of
$z_i(t)$ and $\varphi_i$ instead of
$\varphi_i(t)$. We can write the LHS of the previous inequality in the
form
\[
\sum\limits_i p_i f(y_i)
\]
where $p_i = \dfrac{z_i}{\sum\limits_k z_k}$, $\quad y_i =
-\sum\limits_{j \in V_i} \varphi_j - \ln \dfrac{1}{\varphi_i}
\quad$ and $\quad f(z) = e^z$. Function $f$ is convex and
$\sum\limits_i p_i =1$, hence $\sum\limits_i p_i f(y_i) \ge
f(\sum\limits_i p_i y_i)$ and
\begin{equation} \label{convexity}
\frac{\sum\limits_i z_i \varphi_i \exp\left\{-\sum\limits_{j \in
V_i} \varphi_j\right\}}{\sum\limits_k z_k} \ge
\exp\left\{-\sum\limits_i \frac{z_i}{\sum\limits_k z_k}
\sum\limits_{j \in V_i} \varphi_j - \sum\limits_i
\frac{z_i}{\sum\limits_k z_k} \ln \frac{1}{\varphi_i}\right\}.
\end{equation}
Now consider $\sum\limits_i \dfrac{z_i}{\sum\limits_k z_k}
\sum\limits_{j \in V_i} \varphi_j$ and $\sum\limits_i
\dfrac{z_i}{\sum\limits_k z_k} \ln \dfrac{1}{\varphi_i}$
separately:
\begin{equation} \label{first_summand}
\sum\limits_i \frac{z_i}{\sum\limits_k z_k} \sum\limits_{j \in
V_i} \varphi_j = \frac{\sum\limits_i z_i \sum\limits_{j \in V_i}
\varphi_j}{\sum\limits_k z_k} = \frac{\sum\limits_j \varphi_j
\sum\limits_{i \in V_j} z_i}{\sum\limits_k z_k} =
\frac{\sum\limits_j z_j}{\sum\limits_k z_k} = 1
\end{equation}
(we use the identity $\varphi_j \sum\limits_{i \in V_j} z_i = z_j$
and the symmetry of the neighbourhood relation: $j \in V_i$ iff $i
\in V_j$.)

Note that the logarithmic function is concave. Hence
\begin{equation} \label{second_summand}
\sum\limits_i \frac{z_i}{\sum\limits_k z_k} \ln
\frac{1}{\varphi_i} \le \ln\left(\sum\limits_i
\frac{z_i}{\sum\limits_k z_k} \frac{1}{\varphi_i}\right) =
\ln\left(\frac{\sum\limits_i \frac{z_i}{\varphi_i}}{\sum\limits_k
z_k}\right) = \ln\left(\frac{\sum\limits_i \sum\limits_{j \in V_i}
z_j}{\sum\limits_k z_k}\right) = \ln V.
\end{equation}

Inequality (\ref{inequality}) follows now from (\ref{convexity}),
(\ref{first_summand}) and (\ref{second_summand}).

\end{proof}

\begin{proof}
[Proof of Theorem \ref{th:fluidstable}.]

Corollary \ref{cor:liminf}  implies that $z_i(t)
> 0$ for all $i=1, \ldots, K$ and all $0<t<\inf\{u: z(u)=0\}$. Then we can use
Lemma \ref{le:lyapunov}. Note also that $\sum\limits_i x_i \ge
\sqrt{\sum\limits_i x_i^2}$ for any positive values of $\{x_i\}$.
Hence, Lemma \ref{le:lyapunov} implies that
\[
\left(\sum\limits_i z_i^2(t)\right)' \le - \varepsilon
\sqrt{\sum_i z^2_i(t)}.
\]
Then
\[
\left(\sqrt{\sum\limits_i z_i^2(t)}\right)' \le - \varepsilon/2,
\]
and the result follows.

\end{proof}

\begin{remark} \label{rem:Stolyar}

By applying methods of \cite{Stolyar2}, we can get a similar (but
less explicit) stability result in a more general situation.
Assume now that the system may be asymmetric, i.e. that values of
$\lambda_i$ may differ for different $i$ and the graph
$\mathcal{G}$ may be irregular.

Let $$M = \{\mu: \mu_i = p_ie^{- \sum_{j \in V_i} p_j}, i =
1,\ldots, K, \quad \text{for some} \quad \overline p = (p_1,
\ldots, p_K) \quad \text{with} \quad p_i \ge 0\}.$$

One can show that the vector $(\varphi_1, \ldots, \varphi_K)$ with
$\varphi_i = \dfrac{z_i}{\sum_{j \in V_i} z_j}$ maximises the
function $\sum\limits_{i=1}^K z_i \ln \mu_i$ over all vectors $\mu
\in M$. Based on that, one can obtain the following.


\begin{theorem} \label{th:stolyar}
Assume that there exists a vector $\mu \in M$ such that $\lambda <
\mu$ component-wise. Then the Markov chain $W_n$ is positive
recurrent.
\end{theorem}

A proof of Theorem \ref{th:stolyar} follows the lines of the proof of
Theorem 4 in \cite{Stolyar2}.

\end{remark}

\section{Rate of convergence} \label{sec:rates}

In this section, we again consider the symmetric case. We will
obtain power rates of convergence of $W_n$ to its stationary
distribution in the total variation norm. We expect that one can
similarly prove the geometric ergodicity of the underlying Markov
chain given the light-tailedness of the distribution of the
increments $\{A(n)\}$.

Define the total variation distance between distributions $\pi_1$
and $\pi_2$ by
\[
||\pi_1(\cdotp) - \pi_2(\cdotp)|| = \sup\limits_{|g| \le 1}
\left|\int g(y) \pi_1(dy) - \int g(y) \pi_2(dy)\right|.
\]

\begin{theorem} \label{rates}
Assume that $\lambda < e^{-1}/V$ and $\mathbf E A_i(n)^{p+1} <
\infty$ for some $p \ge 1$ and for all $i=1, \ldots, K$ and $n$.
Assume also that $\mathsf P(A_1(1)=0, A_2(1)=0,...,A_K(1)=0) > 0$.
Then
\[
\lim\limits_{n \to \infty}  n^{p} ||\mathsf P^n(x,\cdotp) -
\pi(\cdotp)|| = 0, \quad x \in {\mathbb N}^K,
\]
where $\mathsf P^n(x, \cdotp)$ is the distribution of $W^x(n)$ and
$\pi(\cdotp)$ is the stationary distribution of $W$.
\end{theorem}

\begin{proof}
[Proof of Theorem \ref{rates}]

The proof of Theorem \ref{rates} is based on the following lemma
which is an analogue of Proposition 5.3 of Dai and Meyn \cite{DaiMeyn}.

\begin{lemma} \label{prop53DM}
Assume that the conditions of Theorem \ref{rates} are satisfied.
Then, for some constants $c < \infty$, $\delta > 0$ and a finite
set $C$,
\[
\mathbf E \left(\sum_{n=0}^{\tau_C(\delta)} |W^x(n)|^p\right) \le
c |x|^{p+1}
\]
for any $x \in \mathbb N^K$, where $\tau_C(\delta) = \min(n \ge
\delta: W(n) \in C)$.
\end{lemma}

\begin{proof}
[Proof of Lemma \ref{prop53DM}]

The proof of Lemma \ref{prop53DM} follows the lines of the proof
of Proposition 5.3 of \cite{DaiMeyn}.

It follows from Theorem \ref{th:fluidstable} that there exists
$t_0$ such that
\[
\lim\limits_{|x| \to \infty} \frac{W^x(|x|t_0)}{|x|} = 0 \quad
\text{a.s.}
\]
Note also that the family of random variables
$\left\{\dfrac{\left|W^x(|x|t_0)\right|^{p+1}}{|x|^{p+1}}\right\}$
is uniformly integrable, since
\[
\frac{\left|W^x(|x|t_0)\right|^{p+1}}{|x|^{p+1}} \le
\frac{\left(\sum_{m=0}^{|x|t_0} \sum_{i=1}^K
A_i(m)\right)^{p+1}}{|x|^{p+1}} \le t_0^{p+1}
\frac{\sum_{m=0}^{|x|t_0} \left(\sum_{i=1}^K A_i(m)
\right)^{p+1}}{|x| t_0}
\]
and the family $\left\{\dfrac{\sum_{m=0}^{|x|t_0}
\left(\sum_{i=1}^K A_i(m) \right)^{p+1}}{|x| t_0} \right\}$ is
uniformly integrable. This is guaranteed by the existence of
$\mathbf E A_i(m)^{p+1}$ for all $i = 1, \ldots, K$ and for all
$m$. Hence,
\[
\lim\limits_{|x| \to \infty} \frac{\mathbf E
\left[|W^x(|x|t_0)|^{p+1}\right]}{|x|^{p+1}} = 0.
\]
Choose $L$ such that
\begin{equation} \label{bound_expect}
\mathbf E \left[|W^x(|x|t_0)|^{p+1}\right] \le \frac12 |x|^{p+1}
\end{equation}
for $|x| \ge L$. Define, as in the proof of Proposition 5.3 of
\cite{DaiMeyn}, the sequence of stopping times $\sigma_0=0,
\sigma_1=t(x)$, and $\sigma_{k+1}=\sigma_k+\theta_{\sigma_k}
\sigma_1, k \ge 1$, where $t(x) = t_0 \max(L, |x|)$ and
$\theta$ is
shift operator on the sample space. We assume that $t_0$ is an
integer. The stochastic process $\hat W_k = W(\sigma_k)$ is a
Markov chain with the transition kernel
\[
\hat P(x, A) = \mathsf P(W^x(t(x)) \in A).
\]
Now (\ref{bound_expect}) implies that
\[
\mathbf E \left\{|\hat W_1|^{p+1} - |\hat W_0|^{p+1} | \hat W_0 =
x \right\} \le - \frac12 |x|^{p+1} + b {\bf I}_C(x),
\]
for the set $C = \{x: |x| \le L\}$ and for some constant $b$. The
Comparison Theorem (Meyn and Tweedie \cite{MeynTweedie}, p. 337) yields that
\begin{equation} \label{comptheorem}
\mathbf E\left[\sum_{n=0}^{k_*-1} |W^x(\sigma_k)|^{p+1}\right] =
\mathbf E\left[\sum_{n=0}^{k_*-1} |\hat W(k)|^{p+1}\right] \le 2
\left\{|x|^{p+1} + b {\bf I}_C(x) \right\}
\end{equation}
where $k_* = \min(k \ge 1: \hat W(k) \in C\}$. To prove Lemma
\ref{prop53DM}, we first show that for some constant $c_0$
\begin{equation} \label{prop53_1}
\mathbf E \left[\sum\limits_{n=\sigma_k}^{\sigma_{k+1}}
|W^x(n)|^p| \mathcal F_{\sigma_k}\right] \le
c_0W^x(\sigma_k)^{p+1}
\end{equation}
which by the strong Markov property amounts to
\[
\mathbf E \sum\limits_{n=0}^{t(x)} |W^x(n)|^p \le c_0 |x|^{p+1}.
\]
This follows from the fact that
\[
\sum\limits_{n=0}^{t(x)} |W^x(n)|^p \le \sum\limits_{n=0}^{t(x)}
\left(\sum_{m=0}^n \sum_{i=1}^{K} A_i(m)\right)^p \le
\sum\limits_{n=0}^{t(x)} \left(\sum_{m=0}^{t(x)} \sum_{i=1}^K
A_i(m)\right)^p
\]
a.s., and from our assumption that $\mathbf E A_i(m) < \infty$ for
all $i=1, \ldots, K$ and for all $m$. Substituting
(\ref{prop53_1}) into (\ref{comptheorem}), we have
\[
\mathbf E \left[\sum_{k=0}^{\infty} \mathbf
E\left[\sum_{n=\sigma_k}^{\sigma_{k+1}} |W^x(n)|^p| \mathcal
F_{\sigma_k}\right] {\bf I}{k < k_*}\right] \le c |x|^{p+1}.
\]
By Fubini theorem and the smoothing property of the
conditional expectation, the LHS is precisely $\mathbf E
\left[\sum_{n=0}^{\sigma_{k_*}} (1+|W^x(n)|^p)\right]$. The
proposition now follows from the fact that $\sigma_{k_*} \ge
\tau_C(t_0 L)$.

\end{proof}

We now apply Proposition 5.4 of \cite{DaiMeyn} with $t=1$. In our
 case, it gives the following bound:
\begin{equation} \label{prop54}
\mathbf E \left\{V(W(1)) - V(W(0)) | W(0) = x \right\} \le - f(x)
+ \kappa
\end{equation}
with $V(x) = \mathbf E \left(\sum_{n=0}^{\tau_C(\delta)}
|W^x(n)|^p \right)$ and $f(x) = |x|^p$.


Further, Lemma \ref{prop53DM} implies that $V(x) \le c |x|^{p+1}$.
Now (\ref{prop54}) yields that
\[
\mathbf E \left\{V(W(1)) - V(W(0)) | W(0) = x \right\} \le
V(x)^{\frac{p}{p+1}} + b {\bf I}_C
\]
for the set $C = \{x: |x| \le L \}$ and for some constant $b$. The
result now follows from Theorem 2.5 of Douc et al. \cite{DFMS}.

\end{proof}

\section{Local stability of fluid limits on the positive orthant}
\label{sec:local}

In this Section we investigate the behaviour of the solution to
the differential equation (\ref{eq:diffeq}). In particular, we
show that if the input rate $\lambda$ is sufficiently small, then
the diagonal is locally unstable.

\subsection{Orbits of the fluid limits}

Recall that $H = \{x \in \mathbb R^K: x_i > 0 \quad \text{for all}
\quad i=1, \ldots, K\}$ is the interior of the positive orthant
and $\ind = (1,...,1)$. For $z(t)$ in $H$, the differential
equation (\ref{eq:diffeq}) may be restated in a closed form as
\begin{equation}
\label{eq:eqdiff}  z' (t) = F (\phi (z(t))),
\end{equation}
with $F_i(x) = \lambda - x_i e^{-\sum_{j \in V_i} x_j}$.  Let
$\Delta = \{x \in H : x_1 = x_2 ... = x_K\}$ be the main diagonal
and $C_u = \{x \in H : | x/|x| - \ind / K | \leq u \}$, $u
> 0$, $C_u$ a cone with direction $\Delta$.
We note that the main diagonal is an orbit of the 
equation $F  (\phi (c \ind)) = (\lambda - e^{-1}/V) \ind$.

Let $A$ be the adjacency matrix of $\mathcal{G}$ and $\{\nu_1,..., \nu_K\}$
its eigenvalues with $\nu_i \leq \nu_{i+1}$. The spectral gap
$\gamma$ is defined by:
$$
\gamma = \min_{i < K} (\nu_K- \nu_i) = \nu_K - \nu_{K-1}.
$$
Note that since $\mathcal{G}$ is $(V-1)$-regular, $\nu_K = V$. The main
result of this section is the following.

\begin{theorem}
\label{th:eqdiff}  If $\lambda > \frac{e^{-1}}{V} ( 1 -
\frac{\gamma^2}{ V^2})$, then there exists $u
> 0$ such that, for all solutions $t \to z(t)$ of equation
(\ref{eq:eqdiff}) with the initial condition in $C_u$,
$$
\lim_{t \to +\infty} \phi ( z(t) ) = \ind / V.
$$
If $\lambda < \frac{e^{-1}}{V} (1-\frac{\gamma^2}{ V^2})$, then
the diagonal is locally unstable.
\end{theorem}

Theorem \ref{th:eqdiff}  will be proved in the next Subsection.

\begin{corollary}
\label{cor:eqdiff} Assume that $\lambda > \frac{e^{-1}}{V}
(1-\frac{\gamma^2}{ V^2})$ and that $z(t)$ is a solution of
(\ref{eq:eqdiff}). There exists $u>0$ such that if $z(0) \in
C_u$, then
\begin{enumerate}
\item[(i)] if $\lambda < e^{-1}/V$, then $z(c) = 0$ for some $c>0$.
\item[(ii)] if $\lambda > e^{-1}/V$, then  $z(t)
\sim (\lambda- e^{-1}/V)t.$
\end{enumerate}
\end{corollary}

\begin{proof}
[Proof of Corollary \ref{cor:eqdiff}]

Let $z(t)$ be the maximal solution with given initial value
$z(0) \in H$. From Theorem
\ref{th:eqdiff}, $ \lim \phi(z(t)) = \ind/V $. Since $F$ is
continuous in a neighbourhood of $\ind/V$, $\lim_{t \to +\infty}
 z'(t) = (\lambda- e^{-1}/V) \ind $. If $\lambda \neq
e^{-1}/V$, then the latter implies that, as $t$ tends to infinity,
\begin{equation}
\label{eq:sim} z(t) \sim  (\lambda- e^{-1}/V) t \ind .
\end{equation}

Then the second statement of the corollary follows. Suppose
now that $\lambda- e^{-1}/V < 0$. Then, from (\ref{eq:sim}), $z(t)$
leaves $H$ in finite time. Lemma \ref{le:trajH} implies in turn
that there exists $c>0$ such that $z(c) = 0$. So, the first assertion
of Corollary \ref{cor:eqdiff} is also proved. \end{proof}


\subsection{Proof of Theorem \ref{th:eqdiff}}

The proof of Theorem \ref{th:eqdiff} is an application of the
stability theory of differential equations. It will be given in
the series of technical lemmas.

\subsubsection{Spectral analysis.}
We need to consider  the eigenvalues of $D (F\circ \phi)(x)$ for
$x \in \Delta$, where $D f (x)$ is the differential of $f$ at $x$.
Here $F\circ \phi$ is homogeneous of order $0$: for all $c>0$,
$F(\phi(cx)) = F (\phi (x))$. Hence,
$$ D (F\circ \phi)(c\ind) = c^{-1} D(F\circ \phi)(\ind).$$ Since
$\Delta$ is an orbit of equation (\ref{eq:eqdiff}), $\ind$ is an
eigenvector of $D (F\circ \phi)(\ind)$ associated with the
eigenvalue $0$.

\begin{lemma}
\label{le:spectrum} The eigenvalues of $D(F\circ \phi)(\ind)$ are
$(0,\eta_1,\cdots, \eta_{K-1})$ with $\eta_i =  - \frac {e^{-1}} {
V ^3 } (V-\nu_{K-i})^2$ . In particular, $\eta_i < 0$ for all $i
\ge 1$.
\end{lemma}

\begin{proof}
[Proof of Lemma \ref{le:spectrum}]

A direct computation leads to: $$ (D (F\circ \phi) (\ind))_{ij} =
\left\{
\begin{array}{ccl}
-\frac{e^{-1}(V-1)}{V^2} &\hbox{if} & j = i, \\
  \frac{e^{-1}}{V^3}|V_i \cup V_j| & \hbox{if}& j \in V_i\setminus \{i\}, \\
 - \frac{e^{-1}}{V^3}|V_i \cap V_j| &  \hbox{if}& j  \not\in V_i.
\end{array}  \right.
$$

Then $D (F\circ \phi) (\ind).\ind = 0$. To show this, let $M = -e V^3 D
(F\circ \phi) (\ind)$. Using the equality $|V_i \cup V_j| = 2V -
|V_i \cap V_j|$, we deduce that: $$ (M{\bf 1})_i = V(V - 1) -
2V(V-1) +\sum_{j\neq i} |V_i \cap V_j| = \sum_{j=1}^K |V_i \cap
V_j| - V^2 = 0. $$

Let $E$ denote the identity matrix and $A$ the adjacency matrix of
$\mathcal G$. Since ($A^2)_{ij} = |V_i \cap V_j|$, we have the
following decomposition: $$ M = V^2 E - 2V A + A^2 = (A-VE)^2. $$

The matrix $A$ is irreducible since $\mathcal G $ is connected.
Thus $(A-VE)$ is an ML-Matrix (refer to Seneta \cite{seneta}). In
the graph theory, this matrix is called the Laplacian matrix of
$\mathcal{G}$. From Corollary 1 of Theorem 1 in Seneta \cite{seneta}, the
spectral radius of $A$ is $V$. Theorem 2.6 (d) of \cite{seneta}
implies that $\dim Ker (A-VE) = 1$ and that all non-zero eigenvalues of
$(A-E)^2$ are positive reals (recall that the
spectrum of $A$ is real). \end{proof}

\subsubsection{Orbit of $\psi \circ z$.}

We define: $$\Sigma =  \{ x \in H \,: \, \sum_{i=1} ^ K x_i = 1 \}
= H \cap \langle \ind ,\cdot \rangle ^{-1} (\{1\}) = \psi (H),$$
where $\psi(x) = x /|x|$. $\Sigma$ is clearly a
$C^{\infty}$-convex manifold of codimension $1$. Introduce the
following differential equation on $\Sigma$:

\begin{equation}
\label{eq:eqdiffS} y'  = D\psi (y) F(\phi(y)) = \alpha(y)
\end{equation}
with an initial condition $y(0)$ in $\Sigma$. Here $\alpha$ is a
$C^\infty(\Sigma)$ function and $\alpha(y) \in T_y (\Sigma)$ the
tangent space of $\Sigma$ at $y$. The next step is to compare the
orbits of equations (\ref{eq:eqdiffS}) and (\ref{eq:eqdiff}). The
next lemma asserts that the orbits of the solution of equation
$y' =
\alpha(y)$ and $\psi \circ x$ are indeed equal (here $t \mapsto z(t)$
is a solution of equation (\ref{eq:eqdiff})).

\begin{lemma}
\label{le:scaling} Let $z(0)$ belong to $H$ and let $z(t)$ be the
maximal solution of equation (\ref{eq:eqdiff}). Let $y(t)$ be the
maximal solution of $ y' = G(y)$, with the initial condition
$y(0) = \psi (z(0))$. Then $y(t)$ is defined on $\mathbb R_+$, and
there exists an increasing
 continuous bijective function $\mu: \mathbb R_+ \to \mathbb R_+$ such that
$$ y \circ \mu  = \psi \circ z $$.
\end{lemma}

\begin{proof}
[Proof of Lemma \ref{le:scaling}]

For any initial condition in $H$, we have $F(\phi(z(t))) \leq \lambda
\ind $. This is clearly true if  $z(t) \in H$. If $z(t) \not\in
H$, then $z(t) \in \Delta \cap H^c$, by Lemma \ref{le:trajH}.
Thus $F(z(t)) = (\lambda - 1/V e^{-1} ) \ind \leq \lambda
\ind$. It follows that $|z(t)| = \sum_{j =1}^K  z_j (t) \leq K
\lambda t + \sum_{j =1}^K  z_j (0)$.

Suppose now that $z(t) \in H$ for all $t$. Then the integral
$\int_0
^{\infty} \left(\sum_{j =1} ^ K z_j (s)\right)^{-1} ds$ diverges. By the
intermediate value theorem, we deduce that there exists an
increasing continuous function $\nu$ such that
\begin{equation}
\label{eq:nu} \int_0 ^ {\nu(t)} \frac {ds}{\sum_{j =1} ^ K z_j
(s)} = t \quad \text{for all} \quad t \geq 0.
\end{equation}
In
particular,
\begin{equation*}
 \nu'(t) = \sum_{j =1} ^ K z_j (\nu(t)).
\end{equation*}

Let $w = \psi \circ z \circ \nu$, with $w(0) = \psi (w(0))= y(0)$.
We have

\begin{eqnarray*}
 w' (t) & = &  \nu'(t) \frac{d}{ds} \psi( z(s)) \Bigm|_{s =
\nu(t)}
\\ & = & \left(\sum_{j =1} ^ K z_j (\nu(t)) \right) D\psi (z(\nu(t)))
F(w(t)).
\end{eqnarray*}

The function $\psi$ is homogeneous of order $0$ and thus $D\psi(c
z) = c^{-1} D\psi(z)$, for all $c>0$. Then

\begin{eqnarray*}
w' (t) & = &  D\psi \left(\frac{z(\nu(t))}{\sum_{j =1} ^ n z_j
(\nu(t))} \right) F(w(t))\\ & = & G(w(t)).
\end{eqnarray*}

The solution of the differential equation is unique, therefore $w
(t) = y(t)$. The lemma is proved with $\mu = \nu ^ {-1}$.


\end{proof}



\subsubsection{Local stability of $\psi \circ z$.}

Clearly, $y_0 = \ind/K$ is an equilibrium point of equation
(\ref{eq:eqdiffS}). In the next lemma we prove that this
equilibrium is locally stable.

\begin{lemma}
\label{le:locstableS} If $\lambda > \frac{e^{-1}}{V}
(1-\frac{\gamma^2}{ V^2})$, there exists $u > 0$ such that, for
all solutions $t \mapsto y(t)$ of equation (\ref{eq:eqdiffS}) with
$|y(0) - y_0 | < u$,
$$
\lim_{t \to +\infty} \sup_{y(0) \in \Sigma : | y(0) - y_0| < u } | y(t) - \ind / V | = 0.
$$
\end{lemma}

\begin{proof}
[Proof of Lemma \ref{le:locstableS}]

We denote by $D \alpha (y) |_{T_y (\Sigma)}$ the differential of
$\alpha$ at $y$ restricted to the $K-1$ dimensional subspace $T_y
(\Sigma)$. It is known that if all the eigenvalues of $D \alpha
(y_0)|_{T_y (\Sigma)}$ have a negative real part, the local
stability follows (see, e.g., Coddington and Levinson
\cite{coddington}). Let $D^2 \psi(y) (\cdot,\cdot)$ denote the
second differential of $\psi$ at $y$, seen as a bilinear mapping.
We have
\begin{equation}
\label{eq:Dalpha} D \alpha (y) = D^2 \psi (y) ( F (\phi (y)),
\cdot ) + D \psi(y) D (F \circ \phi) (y).
\end{equation}
The first term in (\ref{eq:Dalpha}) is a matrix and its entry
$(i,j)$ is equal to
$$
\sum_{k=1} ^K \frac{\partial^2 \psi (y) _i }{\partial y_j
\partial y_k} F (\phi (y)) _k.
$$
Clearly, $F (\phi (y_0)) = (\lambda - e^{-1}/V)\ind$.
Then a straightforward computation gives
$$
D^2 \psi (y_0) ( F (\phi (y_0)), \cdot ) = (\lambda - e^{-1}/V) (J
- K E ),
$$
where $E$ is the identity matrix and $J$ is the matrix with all
its entries equal to $1$. We also have $D\psi(y_0) =E - J/K$.
Finally, equation (\ref{eq:Dalpha}) can be rewritten as
$$
D \alpha (y_0) =  ( E - J / K  )  \Bigm( D (F \circ \phi) (y_0) -
(\lambda - e^{-1}/V) E \Bigm).
$$
The matrix $(E - J/K)$ commutes with all symmetric matrices and
has two eigenvalues: eigenvalue $1$ (with multiplicity $K-1$) and
eigenvalue $0$ (with
multiplicity $1$, associated to the eigenvector $\ind$). By Lemma
\ref{le:spectrum}, the eigenvalues of $D (F \circ \phi) (y_0) -
(\lambda - e^{-1}/V) E$ are
$$ \mu_i = - e^{-1} ( V - \nu_{K-i})^2  /V^3  -   \lambda +  e^{-1}/V
\quad \mbox{for} \quad 0\le i\le K-1.
 $$
The eigenvector associated to $\mu_0 = \lambda -  e^{-1}/V $ is
$\ind$. Thus we have proved that $\lambda - e^{-1}/V$ is an
eigenvalue of multiplicity $1$ for $D \alpha (y_0)$ and that the
other eigenvalues are $(\mu_i)_{i \geq 1}$. These eigenvalues have
negative real parts if and only if $\mu_1 = -e^{-1}\gamma^2  / V^3
- \lambda + e^{-1}/V < 0$, which is equivalent to
$\lambda > e^{-1} ( 1 -\gamma^2 / V^2 ) / V$. The vector space generated by the
associated eigenvectors is precisely the tangent hyperplane
$T_{y_0}( \Sigma) = \ind^\bot$, the hyperplane orthogonal to
$\ind$. \end{proof}

Now we can prove Theorem \ref{th:eqdiff}. Let $|z(0)| \in C_u$ and
 $y(0) = z(0)/|z(0)|$. Then, by Lemmas \ref{le:scaling} and
\ref{le:locstableS}, $$ \lim_{t \to + \infty} \psi(z(t)) = \lim_{t
\to +\infty} y ( \mu (t)) = \ind/K. $$

In particular, $\phi(z(t))$ tends towards $\ind/V$ as $t$ tends to
infinity.



\section{Absence of attraction to the diagonal in one particular case} \label{behaviour}

It has already been pointed out in Section 1.3 and in the
previous Section that if $\lambda$ is too small,
the main diagonal may not be locally stable. In this Section, we present an example
of a graph with locally stable sets of parameters
which do not belong to the main diagonal if $\lambda$ is sufficiently small.

Consider a graph $\mathcal{G}$ with $4$ vertices placed on a circle.
Number the vertices $1,2,3,4$ clockwise and assume
that each vertex is linked with its $2$ neighbours
(so that, for example, vertex $1$ has links with $2$ and $4$). In
this case, $K=4$ and $V=3$.


For the fluid limits associated with this graph, consider equation
(\ref{eq:eqdiffS}). It can be rewritten in the form
$${y}'_i(t) = \left(\lambda - \varphi_i(t) e^{-\sum_{j \in V_i}
\varphi_j(t)}\right) - y_i \sum_{k=1}^K \left(\lambda -
\varphi_k(t) e^{-\sum_{j \in V_k} \varphi_j(t)}\right), \quad
i=1,\ldots,K.$$

We are interested in the so-called stable points of the latter
system of differential equations, i.e. points for which the
RHS's of all the equations above are identically $0$. So if $y(0)$ is
such a point, $y(t)$ stays at this point for all $t \ge 0$.
Clearly, one stable point is $(1/K, \ldots, 1/K)$, which
corresponds to the diagonal. However, if $\lambda <
\dfrac{e^{-1}}{V}\left(1-\dfrac{\gamma^2}{V^2}\right)$ ($=5/27
e^{-1}$ in our case), then there exist other stable points.

Take $y_1(0) = y_2(0)$ and $y_3(0) = y_4(0)$. Since $y_1(0) +
y_2(0) + y_3(0) + y_4(0) = 1$, the equality $y_3(0) =
(1-2y_1(0))/2$ holds, and the system of differential equations at
time $t=0$ reduces to just one (i.e. any) of them. One can show
that, for any $\lambda < 5/27 e^{-1}$, the RHS of this equation
equals $0$ at three different points: at
$\overline{y^{(1)}}=(1/4,1/4,1/4,1/4)$ and at two others, say
$\overline{y^{(2)}}$ and $\overline{y^{(3)}}$. One can find
approximate values of these points numerically. For instance, if
$\lambda = 0.001$, then $\overline{y^{(2)}} \approx (0.01, 0.01,
0.49, 0.49) $ and $\overline{y^{(3)}} \approx (0.49, 0.49, 0.01,
0.01)$. Numerical results also show that these points are locally
stable.

\section{Extensions of the model} \label{sec:extensions}

\subsection{Random neighbourhood}

In this Subsection we consider a possible extension of our model.
Assume there is a fixed number of points $1, \ldots, K$ and a set
of undirected graphs $\left\{\mathcal{G}^j\right\}_{j=1}^L$ each
having points $1, \ldots, K$ as its vertices. Assume that at each
time $n$ the neighbourhood relations are given by the graph
$\mathcal{G}^{\eta_n}$ where $\eta_n$ are independent identically
distributed random variables taking the value $j$ with probability
$p_j$. The need to consider such a variability of neigbourhood
relations may be justified by, for instance, the fact that a
change of environment conditions may lead to a change of the
radius and/or direction of interactions.

Denote by $\mathcal{V}_i^j$ the neighbourhood of the point $i$ in
the graph $\mathcal{G}^j$ and by $V_i^j$ its cardinality. We assume
again that the system is regular in the sense that $\mathbf E V_i^{\eta_1}
= V$ for all $i$.

Following the proof of Theorem \ref{proofderivative}, one can show
that the fluid limits of the model described above satisfy the
following differential equation
\begin{equation*}
\label{eq:diffeq1} z_i'(t) = \lambda - \sum\limits_{k=1}^L p_k
\varphi_i^k(t) e^{- \sum\limits_{j \in V_i^k} \varphi_j^k(t)}
\end{equation*}
where $\varphi_i^k(t)$ are defined in a natural way. Using the
same methods as those used in the proof of Theorem
\ref{th:fluidstable}, it can be shown that the system with random
neighbourhood is stable if $\lambda < \dfrac{e^{-1}}{V}$.

\subsection{Non-regular graphs with space-inhomogeneous input}

Although Remark \ref{rem:Stolyar} provides sufficient conditions
for stability in this case, the conditions are not easy to verify. Here we
give some other conditions that are also sufficient for the
stability of the system. Assume now that $\mathbf E V_i^{\eta_1} =
V_i$, and $V_i$ are not necessarily equal. Put $V = \max\limits_i
V_i$. Assume also that $\mathbf E \xi_i^n = \lambda_i$, so that
the input is ``space-inhomogeneous''. Put $\lambda = \max\limits_i
\lambda_i$. Clearly, all the results concerning fluid limits also
hold in this case, and it is easy to see that one can prove the
following result.

\begin{theorem}
The system described above is stable provided $\lambda <
\dfrac{e^{-1}}{V}$.
\end{theorem}

\acks
The authors are grateful to F.~Baccelli and D.~Denisov for useful
discussions on different parts of this work, to A.~Stolyar for
bringing his paper \cite{Stolyar2} to their attention, and to
A.~Richards for careful reading and important comments.

The authors would also like to thank Framework 7 EURO NF project for
the travel support. The second author acknowledges the travel support
from EPSRC grant No.~R58765/01.

\end{document}